\documentclass[12pt]{amsart}

\usepackage{graphicx,amscd}
\input xy
\xyoption{all}

\newsymbol\pp 1275
\newsymbol\twoheadrightarrow 1310

    \newcommand{\A}{{\mathbb A}}
\newcommand{\Sym}{\operatorname{Sym}}
\newcommand{\Grass}{\operatorname{Grass}}
\newcommand{\Hom}{\operatorname{Hom}}
\newcommand{\Ext}{\operatorname{Ext}}
\newcommand{\Rep}{\operatorname{Rep}}

\newcommand{\mult}{\operatorname{mult}}

\newcommand{\id}{\operatorname{id}}

\newtheorem{theorem}{Theorem}
\newtheorem{proposition}[theorem]{Proposition}

\newtheorem{lemma}[theorem]{Lemma}
\newtheorem{corollary}[theorem]{Corollary}

\theoremstyle{definition}

\newtheorem{example}[theorem]{Example}

\newcommand{\Z}{{\mathbb Z}}
\newcommand{\GL}{\operatorname{GL}}
\newcommand{\SL}{\operatorname{SL}}

\newcommand{\SI}{\operatorname{SI}}

\title[the number of subrepresentations for general quiver representations]{On the number of subrepresentations\\ of a general quiver representation}
\author{Harm Derksen, Aidan Schofield and Jerzy 
Weyman}
\thanks{The first author is supported by NSF, grant DMS 0349019
and the third author is also supported by NSF, grant DMS 0300064.}
\begin{document}
\maketitle
\section{Introduction}
It is well-known that the intersection multiplicities
of Schubert classes
in the Grassmanian are Littlewood-Richardson
coefficients. For a partition $\lambda$ inside a $r\times (n-r)$
rectangle, let $Y_\lambda$ be the Schubert variety inside the Grassmanian
$\Grass(r,n)$ 
corresponding to $\lambda$
and let $[Y_{\lambda}]$ be its cohomology class.
We have
\begin{equation}\label{eqa}
[Y_\lambda]\cdot [Y_\mu]=\sum_{\nu} c_{\lambda,\mu}^{\nu} [Y_{\nu}]
\end{equation}
where $c_{\lambda,\mu}^{\nu}$ is a Littlewood-Richardson coefficient
and the sum runs over all partitions $\nu$ inside
a $r\times (n-r)$ rectangle. This result can be translated into the language of
quivers as follows. Consider the triple flag quiver 
$$
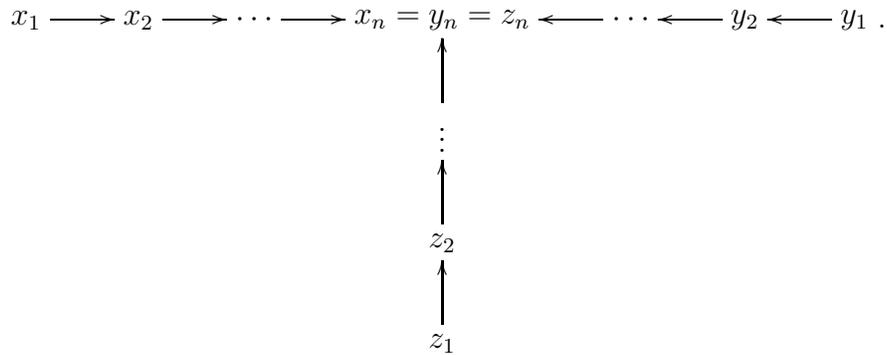
\xymatrix{
x_1\ar[r] & x_2\ar[r] & \cdots\ar[r] &x_n=y_n=z_n& \cdots\ar[l] & y_2\ar[l] &
y_1\ar[l] \\
& & & \vdots\ar[u] & & &\\
& & & z_2\ar[u] & & &\\
& & & z_1\ar[u] & & &
}.
$$
Suppose that $\lambda,\mu,\nu$ are partitions whose Young diagram
fit into a $r\times (n-r)$ rectangle,
such that $|\lambda|+|\mu|+|\nu|=r(n-r)$.
The intersection of the classes $[Y_\lambda]$, $[Y_\mu]$ and
$[Y_\nu]$ is zero-dimensional.
Define dimension vectors $\alpha,\beta$ by
$$
\alpha=\begin{matrix}
1 & 2 & \cdots & n & \cdots & 2 & 1\\
& & & \vdots & & & \\
& & & 2 & & &\\
& & & 1 & & &\end{matrix}
$$
and 
$$
\beta(x_{n-r-\lambda_i+i})=\cdots =\beta(x_{n-r-\lambda_{i+1}+i})=i
$$
$$
\beta(y_{n-r-\mu_i+i})=\cdots =\beta(y_{n-r-\mu_{i+1}+i})=i
$$
$$
\beta(z_{n-r-\nu_i+i})=\cdots =\beta(z_{n-r-\nu_{i+1}+i})=i
$$
for $i=1,2,\dots,r$ with the ad-hoc conventions
$\lambda_{r+1}=\mu_{r+1}=\nu_{r+1}=0$ and 
$\lambda_0=\mu_0=\nu_0=n-r-1$.

The number of $\beta$-dimensional representations
of a {\it general} $\alpha$-dimensional representation is finite
and equal to
the multiplicity of $[Y_{\overline{\nu}}]$ inside
$[Y_{\lambda}]\cdot [Y_{\mu}]$.
Here
$$\overline{\nu}=(n-r-\nu_{r},n-r-\nu_{r-1},\cdots,n-r-\nu_1)$$
 is the complementary partition of $\nu$ inside the $r\times (n-r)$ rectangle.

Following the calculation in \cite{D-W},
the dimension of the space of semi-invariants of weight 
$\langle\beta,\cdot\rangle$
on the space $\Rep(Q,\alpha-\beta)$ of $(\alpha-\beta)$-dimensional
representations is
the Littlewood-Richardson coefficient $c_{\lambda,\mu}^{\overline{\nu}}$,
where $\langle \cdot,\cdot\rangle$ is the Euler form defined in the next section.

In this paper we generalize the connection between Schubert calculus and
Littlewood-Richardson coefficients to  quivers without
oriented cycles and
arbitrary dimension vectors.
For a quiver $Q$ without oriented cycles
and two dimension vectors $\alpha,\beta$ we define $N(\beta,\alpha)$
as the number $\beta$-dimensional subrepresentations of a general 
$\alpha$-dimensional
representation,
and $M(\beta,\alpha)$ as the dimension of the
space of semi-invariants polynomials
of weight $\langle \beta,\cdot\rangle$ on the representation
space of dimension
$\gamma:=\alpha-\beta$. 
\begin{theorem}\label{theo1}
If $\langle \beta,\gamma\rangle=0$, then we have 
$N(\beta,\alpha)=M(\beta,\alpha)$.
\end{theorem}
The proof compares the calculations of $N(\beta,\alpha)$ and $M(\beta,\alpha)$.
The calculation of $N(\beta,\alpha)$ comes from
Intersection Theory and was done by Crawley-Boevey (\cite{CB}). The 
calculation of $M(\beta,\alpha)$
comes from standard calculations in the coordinate ring of a 
representation space involving the Littlewood-Richardson rule.
 In this note we explain why both calculations are the same.
 For this we will make use of (\ref{eqa}).
 
Suppose that a general representation of dimension $\alpha$ has
infinitely many $\beta$-dimensional subrepresentations. In Section~\ref{sec5}
we will see that in that case, the cohomology class of the variety
of $\beta$-dimensional subrepresentations of a $\alpha$-dimensional
representation in general position is given by a formula
whose coefficients can be interpreted as multiplicities
of isotypic components in the coordinate ring of $\Rep(Q,\gamma)$.

Suppose that $\langle \beta,\gamma\rangle=0$. Given a general
representation of dimension $\alpha$, we will construct in Section~\ref{sec6} a
basis of the semi-invariant polynomials on $\Rep(Q,\gamma)$
of weight $\langle\beta,\cdot\rangle$.
 
\section{Basic Notation}\label{sec2}
  A quiver is a pair $Q= (Q_0 ,Q_1)$ where
$Q_0$ is the set of vertices and $Q_1$ is the set of arrows.
Each arrow $a$ has head $ha$ and tail $ta$,
both in $Q_0$:
$$ta{\buildrel a\over\longrightarrow}ha.$$
An oriented cycle is a sequence of arrows $a_1,a_2,\dots,a_r\in Q_1$
such that $ta_i=ha_{i+1}$ for $i=1,2,\dots,r-1$ and $ha_1=ta_r$.
We will assume that $Q$ has no oriented cycles.

We fix an algebraically closed base field $K$. A representation $V$ of $Q$ is a
family of
 finite dimensional
$K$-vector spaces 
$$\lbrace V(x) \mid x\in Q_0\rbrace$$ 
together with a family of $K$-linear maps
$$\{V(a):V(ta)\rightarrow V(ha)\mid a\in Q_1\}.$$ 
The dimension vector of a representation
$V$ is the function ${\underline d(V)}:Q_0\rightarrow \Z$ defined by
$\underline d(V)(x):= \dim V(x)$. The dimension vectors lie in the space
$\Gamma=\Z^{Q_0}$ of integer-valued functions
on $Q_0$. A morphism $\phi :V\rightarrow W$ of two representations
is a collection of $K$-linear maps
$$\{\phi(x): V(x)\rightarrow W(x)\mid x\in Q_0\}$$ 
such that for each $a\in Q_1$ we have
$W(a)\phi(ta)=\phi (ha)V(a)$, i.e., the diagram
$$
\xymatrix{
V(ta)\ar[r]^{V(a)}\ar[d]_{\phi(ta)} & V(ha)\ar[d]^{\phi(ha)}\\
W(ta)\ar[r]_{W(a)} & W(ha)}
$$
commutes.
We denote the vector space of morphisms from $V$ to $W$ by 
$\Hom_Q(V,W)$.


The category of representations of $Q$ is hereditary, i.e., a subobject
of a projective
object is projective. This implies that if $V$ and
$W$ are representations, then $\Ext_Q^i(V,W)=0$ for all $i\geq 2$.
We will write $\Ext_Q(V,W)$ instead of $\Ext_Q^1(V,W)$.

The spaces $\Hom_Q (V, W)$ and $\Ext_Q (V, W)$ can be calculated as the kernel
and cokernel of the following linear map
\begin{equation}\label{eq1a}
d^V_W : \bigoplus_{x\in Q_0} \Hom (V(x),
W(x))\longrightarrow\bigoplus_{a\in Q_1}
\Hom (V(ta), W(ha)).
\end{equation}
where the  map $d^V_W$ restricted to $\Hom (V(x), W(x))$ is equal to
$$
\sum_{\scriptstyle a \atop \scriptstyle ta=x} 
\Hom (\id_{V(x)},W(a))-
\sum_{\scriptstyle a\atop \scriptstyle ha=x} \Hom (V(a),\id_{W(x)}).
$$
In other words,
$$
d^V_W\big(\{\phi(x)\mid x\in Q_0\}\big)= \{W(a)\phi(ta)-\phi(ha)V(a)\mid a\in Q_1\}.
$$

For a dimension vector $\beta$ we denote by 
$$\Rep(Q,\beta)=\bigoplus_{a\in Q_1}\Hom(K^{\beta(ta)},K^{\beta(ha)})$$ 
the vector space of representations
of $Q$ of dimension vector $\beta$. The groups 
$$\GL(Q,\beta):= \prod_{x\in
Q_0} \GL(\beta (x))$$
and its subgroup 
$$\SL(Q,\beta)=\prod_{x\in Q_0} \SL(\beta (x))
$$
act on $\Rep(Q,\beta)$ as follows. If 
$$A=\{A(x)\mid x\in Q_0\}\in \GL(Q,\beta)$$
and 
$$V=\{V(a)\mid a\in Q_1\}\in \Rep(Q,\beta),$$ 
then we define
$$
A\cdot V:=\{A(ha)V(a)A(ta)^{-1}\mid a\in Q_1\}.
$$
The group $\GL(Q,\beta)$ acts on the coordinate ring
$K[\Rep(Q,\beta)]$ as follows. If $f\in K[\Rep(Q,\beta)]$,
$A\in \GL(Q,\beta)$ then
$$
(A\cdot f)(V)=f(A^{-1}\cdot V),\quad V\in \Rep(Q,\beta).
$$

We are interested in the ring
of semi-invariants
$$\SI(Q,\beta) = K[\Rep (Q,\beta)]^{\SL(Q,\beta)}.$$
To each $\sigma\in \Gamma=\Z^{Q_0}$ we can associate
a character of $\GL(\beta)$ defined by
$$
A=\{A(x)\mid x\in Q_0\}\in \GL(\beta)\mapsto 
\prod_{x\in Q_0}\det(A(x))^{\sigma(x)}.
$$
By abuse of notation this character is also denoted by $\sigma$.
The 1-dimensional representation corresponding to this
multiplicative character will be denoted by $\det^\sigma$.
For any two characters $\sigma,\tau\in \Z^{Q_0}$ we have
$\det^{\sigma+\tau}=\det^{\sigma}\otimes \det^{\tau}$.
We can write
$$
{\textstyle \det^\sigma}=\bigotimes_{x\in Q_0}{\textstyle \det_x^{\sigma(x)}}
$$
where $\det_x^k$ is the $1$-dimensional representation of $\GL(\beta(x))$
corresponding to the multiplicative character $A\mapsto \det(A)^k$.

The ring $\SI(Q,\beta)$ has a weight space decomposition
$$\SI(Q,\beta)=\bigoplus_\sigma \SI(Q,\beta)_\sigma$$
where $\sigma$ runs through the characters of $\GL(Q,\beta)$ and
$$\SI(Q,\beta)_\sigma = \lbrace\ f\in K[\Rep(Q,\beta)]\mid g(f)=\sigma
(g)f\ \forall g\in \GL(Q,\beta)\ \rbrace.$$

Let $\alpha ,\beta$ be two elements of $\Gamma$. We define the Euler inner
product
\begin{equation}\label{eq2a}
\langle\alpha ,\beta \rangle = \sum_{x\in Q_0} \alpha (x)\beta (x)-
\sum_{a\in Q_1} \alpha(ta)\beta(ha).
\end{equation}
It follows from (\ref{eq1a}) and (\ref{eq2a}) that
$$\langle\underline d(V),\underline d(W)\rangle= \dim_K
\Hom_Q (V,W)-\dim_K \Ext_Q (V,W).$$

\section{The computation of $N(\beta,\alpha)$}\label{sec3}

If $r,n$ are nonnegative integers with $r\leq n$ then
$\Grass(r,n)$ denotes the Grassmanian of $r$-dimensional subspaces of $K^n$.
Let $\alpha ,\beta,\gamma$ be two dimension vectors such that
$\alpha=\beta+\gamma$. 
We define 
$$\Grass(\beta,\alpha )=\prod_{x\in Q_0}
\Grass(\beta(x), \alpha(x)).$$ 
A point $W=\{W(x)\mid x\in Q_0\}\in\Grass(\beta,\alpha )$
is a collection of subspaces with $\dim W(x)=\beta(x)$ for all $x\in Q_0$.
Consider the incidence variety
$$
Z(Q,\beta,\alpha)=\lbrace\ (V, W)\in \Rep(Q,\alpha )\times
\Grass(\beta,\alpha )\mid
\forall a\in Q_1\ V(a)(W(ta))\subseteq W(ha)\ \rbrace .$$
There are two projections:
$$
\xymatrix{ & Z(Q,\beta,\alpha)\ar[ld]_{q}\ar[rd]^{p} & \\
\Rep(Q,\alpha) & & \Grass(\beta,\alpha)}
$$
The following proposition was proved in~\cite{S1}.
\begin{proposition}\label{prop1}\ 
\begin{enumerate}
\renewcommand{\theenumi}{\alph{enumi}}
\item The first projection  $q:Z(Q,\beta,\alpha)\rightarrow \Rep(Q,\alpha )$ 
is proper.
\item The second projection $p:Z(Q,\beta,\alpha)\rightarrow
\Grass(\beta,\alpha )$ gives
$Z(Q,\beta,\alpha)$ the structure of a vector bundle over
$\Grass(\beta, \alpha )$.
\item  $\dim\ Z(Q,\beta,\alpha)-\dim\ \Rep(Q,\alpha )=\langle \beta
,\gamma \rangle$.
\end{enumerate}
\end{proposition}
We will use a result of Crawley-Boevey
(\cite{CB}) who proved a formula for the cohomology class of the general fiber of $q$ in terms of Intersection Theory.

In order to formulate this result we need some notation. 
For a variety $X$, its associated cohomology ring
will be denoted by ${\mathcal A}^* (X)$ (the Chow ring or the singular
cohomology ring). We recall
that for the Grassmannian $\Grass(r,n)$ the ring ${\mathcal A}^* (\Grass(r,n))$ is
spanned by classes of
Schubert varieties $Y_\lambda$ where $\lambda$ denotes a partition
contained in the rectangle
$((n-r)^{r})$, i.e., $\lambda=(\lambda_1,\dots,\lambda_r)$ 
with $n-r\geq\lambda_1\geq \lambda_2\geq \cdots \geq \lambda_1\geq 0$.
 The sum $|\lambda|:=\lambda_1+\lambda_2+\cdots+\lambda_r$ of the parts of
$\lambda$ is equal to  the codimension of
$Y_\lambda$.
In our setup we have 
$${\mathcal A}^* (\Grass(\beta,\alpha ))=\bigotimes_{x\in Q_0} {\mathcal A}^*
(\Grass(\beta (x),\alpha(x)).$$ 
We
 denote by $[\lambda]_x$ the cohomology class of the Schubert variety
$Y_\lambda$ in the factor ${\mathcal A}^*
(\Grass(\beta (x),\alpha(x)))$. 
We use the convention $[\lambda ]_x =0$ if
$\lambda$ is not contained in the
rectangle $(\gamma (x)^{\beta (x)})$. 

\begin{proposition}[\cite{CB}]\label{prop2}
For general $V\in \Rep(Q,\alpha)$ the cycle of
$q^{-1}(V)$ in
${\mathcal A}^* (\Grass(\beta,\alpha ))$ is equal to
\begin{equation}\label{eq0}
[q^{-1}(V)]=\prod_{a\in Q_1} \big(\sum_{\lambda}
[\lambda]_{ta}[{\overline\lambda}]_{ha} \big)
\end{equation}
where  in the sum $\sum_{\lambda}[\lambda]_{ta}[{\overline\lambda}]_{ha}$,
$\lambda$
runs over all partitions
which fit inside a $\beta(ta)\times \gamma(ha)$
rectangle. In (\ref{eq0}), $\overline{\lambda}$ denotes
the complement of $\lambda$ inside a $\beta(ta)\times \gamma(ha)$ rectangle.
\end{proposition}
We exchange the sum and the product in (\ref{eq0}). 
Clearly the summands in the
formula correspond to
functions ${\underline\lambda} :Q_1 \rightarrow {\mathcal P}$ where $\mathcal P$ is
a set of partitions
and $\lambda(a)$ is contained in a $\beta(ta)\times \gamma(ha)$ rectangle
for all arrows $a$.
 Now (\ref{eq0}) can be rewritten as
\begin{equation}\label{eq1}
[q^{-1}(V)]=\sum_{\underline\lambda :Q_1\rightarrow {\mathcal P}} \prod_{x\in
Q_0} \prod_{{\scriptstyle a\in Q_1\atop \scriptstyle ta=x}}
[{\underline\lambda} (a)]_x  \prod_{{\scriptstyle a\in Q_1\atop
\scriptstyle ha=x}}[{\overline{\underline\lambda }(a)}]_x.
\end{equation}
Suppose that $\langle \beta,\gamma\rangle=0$.
Then we have $\dim\ Z(Q,\beta,\alpha)=\dim\ \Rep(Q,\alpha )$ by
Proposition~\ref{prop1}(c).
This means that
the general fiber of $q$ is finite. 
Such a general fiber is reduced, even in positive characteristic
(see \cite[Corollary 3]{CB}).
Let $N(\beta ,\alpha )$ be the
cardinality of a general
fiber. 
So a general representation $V$ of dimension
$\alpha$ has $N(\beta ,\alpha )$ subrepresentations of dimension $\beta$.
We translate (\ref{eq1}) into the language of Schur functors
\begin{equation}\label{eqstar}
 N(\beta,\alpha)=|q^{-1}(V)|=
 \sum_{\underline\lambda :Q_1\rightarrow {\mathcal P}}
\prod_{x\in Q_0} \mult\big(S^{\gamma(x)^{\beta(x)}},
\bigotimes_{{\scriptstyle a\in Q_1\atop \scriptstyle ta=x}}
S^{{\underline\lambda} (a)}\otimes \bigotimes_{{\scriptstyle a\in Q_1\atop
\scriptstyle ha=x}}S^{{\overline{\underline\lambda}(a)}}\big)
\end{equation}
Here $\mult(S^{\lambda};T)$ denotes the multiplicity of $S^{\lambda}$ in $T$.
\section{The computation of $M(\beta,\alpha)$}\label{sec4}

Let us calculate the dimension of a weight space $\SI(Q,\gamma
)_{\langle\beta ,\cdot \rangle}$.
For now we assume that the base field has characteristic 0. 
This is sufficient, because we will show in the next section
that the number $M(\beta,\alpha)$ does not depend on the (algebraically closed)
base field.
The space $\Rep(Q,\gamma )$ can be identified with
$$\prod_{a\in Q_1} \Hom_K (W(ta), W(ha)).$$
where $W(x)$ is a vector space of dimension $\gamma (x)$ for
all $x\in Q_0$. The coordinate
ring can now be identified
with the symmetric algebra on the dual space
$$K[\Rep(Q,\gamma )]= \bigotimes_{a\in Q_1}\Sym\big(W(ta)\otimes W(ha)^*\big).$$
Here $\Sym$ denotes the symmetric algebra on a vector space.
 By Cauchy's formula we can rewrite this in terms of Schur functors. We use
the exterior power
notation for Schur functors, i.e.,
for a partition $\mu$ we write ${\textstyle \bigwedge^\mu} W:=S^{\mu'}W$
where $\mu'$ is the conjugate partition of $\mu$.
 In particular, 
${\textstyle \bigwedge^{(m)}}$ is the $m$-th
exterior power.
We have
\begin{equation}\label{eq2ab}
K[\Rep(Q,\gamma )]=\bigoplus_{{\underline\lambda}:Q_1\rightarrow{\mathcal
P}}\bigotimes_{a\in Q_1}
\big({\textstyle \bigwedge^{{\underline\lambda}(a)}}W(ta)\otimes
{\textstyle \bigwedge^{{\underline\lambda}(a)}}W(ha)^* \big).
\end{equation}
This can be rewritten as
\begin{equation}\label{eq2}
K[\Rep(Q,\gamma )]=\bigoplus_{{\underline\lambda}:Q_1\rightarrow{\mathcal
P}}\bigotimes_{x\in Q_0}
\big(\bigotimes_{{\scriptstyle a\in Q_1\atop \scriptstyle ta=x}}
{\textstyle \bigwedge^{{\underline\lambda}(a)}}W(x)\big)\otimes 
\big(\bigotimes_{{\scriptstyle a\in
Q_1\atop \scriptstyle ha=x}}{\textstyle \bigwedge^{{\underline\lambda}(a)}}W(x)^* \big).
\end{equation}

Let us  calculate the dimension of the space of semi-invariants of
weight $\langle\beta
,\cdot\rangle$. 

The partition ${\underline\lambda}(a)$ has parts $\le\gamma (ha)$,
because $\dim W(ha)=\gamma(ha)$.
We need additional restriction to match the one in Proposition~\ref{prop2}.
This is provided by the  following lemma.
\begin{lemma}\label{lem1}
 If a summand in (\ref{eq2}) corresponding to the function
$\underline\lambda :Q_1\rightarrow
{\mathcal P}$ contains a nonzero semi-invariant of weight $\langle\beta
,\cdot\rangle$, then for each $a\in
Q_1$ the partition ${\underline\lambda}(a)$ is contained in
$(\gamma (ha))^{\beta (ta)}$.
\end{lemma}
\begin{proof} Let us look at the space of semi-invariants
$\SI(Q,\gamma )_{\langle \beta ,\cdot\rangle}$.
By Theorem~1 of \cite{D-W} (see also \cite{SV}) 
this space is spanned by the semi-invariants $c^V$
defined by the formula
$c^V(W):=\det\ d^V_W$, where $V\in
\Rep(Q,\beta )$ and $d^V_W$ is the differential in (\ref{eq1a}).
 Let us investigate the
contribution of the coefficients of the matrix
$W(a)$ to such a semi-invariant. The only block of $d^V_W$ where this matrix
$W(a)$ occurs is the block
$\Hom (\id_{V(x)},W(a))$:
$$\begin{pmatrix}W(a)&0&\cdots&0\\
0&W(a)&\cdots &0\\
\vdots&\vdots&\ddots&\vdots\\
0&0&\cdots&W(a)\end{pmatrix}$$
with $\beta (ta)$ blocks $W(a)$. Now any multihomogeneous semi-invariant
will come from exhibiting a deteminant
$c^V(W)$ as a polynomial in coefficients of matrices from $V$ and taking a
coefficient of some monomial. Such
a semi-invariant has to be a linear combination of minors of the above
block, multiplied by polynomials
depending on other matrices $W(b)$. But the minors of the above block
matrix are  products
of minors of $W(a)$ with at most $\beta (ta)$ factors in each summand.
By the straightening law (compare for example  \cite{DC-E-P}) we know that the
products of $\beta(ta)$ minors of $W(a)$ are contained in the
space
$$\sum_{\nu} 
{\textstyle \bigwedge^\nu} W(ta)\otimes{\textstyle \bigwedge^\nu} W(ha)^* $$
in $\Sym(W(ta)\otimes W(ha)^* )$, where $\nu$ runs over partitions
with at most $\beta(ta)$ parts.\end{proof}

Let us calculate the dimension of the space of semi-invariants of weight
$\sigma:=\langle\beta
,\cdot\rangle$. We have
\begin{equation}\label{eqsig}
\sigma(x)=\beta(x)-\sum_{\scriptstyle a\in Q_1\atop \scriptstyle ha=x}\beta(ta)
\end{equation}
for all $x\in Q_0$.
By duality, we have the following $\GL(W(ha))$-isomorphism
\begin{equation}\label{eqdual}
{\textstyle \bigwedge^{{\underline\mu}(a)}}W(ha)^* =
{\textstyle \bigwedge^{{\overline{\underline\mu}(a)}}}W(ha)\otimes
{\textstyle\det_{ha}^{-\beta(ta)}}
\end{equation}
where $\overline{\underline{\mu}}(a)$ is the complement
of $\underline{\mu}(a)$ inside an $\beta(ta)\times \gamma(ha)$
rectangle.

From (\ref{eqdual}) and  (\ref{eq2}) follows 
that the dimension $M(\beta,\alpha)$ of $\SI(Q,\gamma
)_{\langle\beta ,\cdot\rangle}$ is equal to
\begin{equation}\label{eqstarstar}
\sum_{\underline\lambda;Q_1\to {\mathcal P}} \prod_{x\in
Q_0} \mult\big({\textstyle\det_x^{\sigma(x)}}; 
\bigotimes_{\scriptstyle a\atop \scriptstyle ta =x}
{\textstyle \bigwedge^{{\underline\lambda}(a)}}W(x)\otimes 
\bigotimes_{\scriptstyle a\atop \scriptstyle ha=x}
({\textstyle \bigwedge^{{\overline{\underline\lambda}(a)}}}W(x)\otimes
{\textstyle\det_x^{-\beta(ta)}})\big )
\end{equation}
The equations (\ref{eqsig}) and (\ref{eqstarstar}) imply
\begin{equation}\label{eqsuper}
M(\beta,\alpha)=\sum_{\underline\lambda;Q_1\to {\mathcal P}} \prod_{x\in
Q_0} \mult\big({\textstyle\det_x^{\beta(x)}}; 
\bigotimes_{\scriptstyle a\atop \scriptstyle ta =x}
{\textstyle \bigwedge^{{\underline\lambda}(a)}}W(x)\otimes 
\bigotimes_{\scriptstyle a\atop \scriptstyle ha=x}
{\textstyle \bigwedge^{{\overline{\underline\lambda}(a)}}}W(x)\big )
\end{equation}
In view of Lemma~\ref{lem1}, we only need to sum over those functions
$\underline{\lambda}:Q_1\to{\mathcal P}$ for which $\underline{\lambda}(a)$ lies
in a $\beta(ta)\times \gamma(ha)$ rectangle.
We write (\ref{eqsuper}) in terms of Schur functors:
\begin{equation}\label{eqstarstar2}
M(\beta,\alpha)=\sum_{\underline\lambda} \prod_{x\in
Q_0} \mult\big({\textstyle \bigwedge^{(\gamma (x)^{\beta (x)})}} ; 
\bigotimes_{\scriptstyle a\in Q_1\atop \scriptstyle ta =x}
{\textstyle \bigwedge^{{\underline\lambda}(a)}}\otimes 
\bigotimes_{\scriptstyle a\in Q_1\atop \scriptstyle ha=x}
{\textstyle \bigwedge^{{\overline{\underline\lambda}(a)}}}\big ).
\end{equation}
\begin{proof}[Proof of Theorem~\ref{theo1}]
The number $N(\beta,\alpha)$ does not depend on the base field,
and neither does $M(\beta,\alpha)$ by Proposition~\ref{propGF3} in the
next section. We may assume that the base field $K$ is algebraically
closed and has characteristic 0.
Note that $S^\lambda=\bigwedge^{\lambda'}$ where $\lambda'$ is the 
conjugate partition. Also, the Littlewood-Richardson coefficients
$c_{\lambda,\mu}^\nu$ and $c_{\lambda',\mu'}^{\nu'}$ are the same.
From these observations it follows that
the righthand-sides in (\ref{eqstar}) and (\ref{eqstarstar}) 
 are
identical.
\end{proof}
\section{Good filtrations}

We work over an algebraically closed field $K$ of arbitrary characteristic.
Suppose that $V$ is an $n$-dimensional vector space. 
We use the convention that a Schur module 
$L_{\lambda'}V$ is denoted by the partition of its highest weight, 
i.e., $\bigwedge^i V=L_{1^i}V$, $S^i V=L_i V$.

Let us recall that {\it a good filtration}
for a rational finite dimensional $\GL(n):= \GL(V)$-module 
$W$ is a filtration
$$0=W_0\subset W_1\subset\ldots\subset W_{s-1}\subset W_s =W$$
for which each factor $W_{i+1}/W_i$ is isomorphic to some
Schur module $L_{\lambda(i)} V$. Such filtration, if it exists,
 may not be unique. The number of factors 
$L_\lambda V$ in any 
good filtration of $W$ does not depend on the choice of a good filtration, 
and will be denoted by $n_\lambda (W)$.
\begin{proposition} \label{propGF1}
If the modules $W_1$ and $W_2$ have good filtrations, 
then the tensor product $W_1\otimes W_2$ has a good filtration.
\end{proposition}
\begin{proof}
 See~\cite[(4.2) Theorem]{Wang} or \cite[Theorem 4.3.1]{Donkin}.\end{proof}
\begin{proposition}\label{propGF2}
Let $W$ be a module with good filtration. Then there exists a good filtration
$$0=W'_0\subset W'_1\subset\ldots\subset W'_{s-1}\subset W'_s =W$$
such that the  submodule of the $SL(V)$-invariants  $W^{\SL(V)}$ in $W$ 
is equal to $W'_t$ for some~$t$. In particular, 
the dimension of $W^{\SL(V)}$ is equal to the 
number of factors $W_i /W_{i-1}$ isomorphic to the trivial representation.
\end{proposition}
\begin{proof}
We use the results of \cite{Donkin} freely. Let us order the highest weights 
$\lambda$ by saying that $\lambda<\mu$ if $\lambda -\mu$ is a sum of 
positive roots for $\SL(V)$. 

Let $W$ be a module with a good filtration. We can assume without loss 
of generality that $W$ is a polynomial homogeneous representation of 
degree $d$. Then by \cite[Proposition 3.2.6]{Donkin} there exists a 
filtration
$$0=W'_0\subset W'_1\subset\ldots\subset W'_{s-1}\subset W'_s =W$$
with $W'_i /W'_{i-1}=L_{\lambda(i)}V$, for which 
$\lambda (1)\le\lambda(2)\le\ldots\le \lambda(s)$. 
We notice that for existence of $\SL(V)$-invariants it is nessesary that 
$d=ne$ for some $e$. Now among all the possible highest weights 
$\lambda$  that correspond to partitions of $d$ with $\leq n$ parts, 
the smallest one is $\lambda = (e^n )$.
Let $t$ be maximal such that $\lambda(1)=\cdots=\lambda(t)=(e^n)$.
All composition
factors of $W'_t$ has  are isomorphic to
the trivial $\SL(V)$-representation $L_{(e^n )}V$.
Hence $\SL(V)$ acts trivially on $W'_t$. 
This proves the first part of the proposition. 

It remains to show that the submodule $W^{SL(V)}$ is equal to $W'_t$.
Let us assume that the opposite is true. Then the
module $W/W'_t$ is a polynomial homogeneous representation of $\GL(V)$ 
which has a good filtration with no factors isomorphic to $L_{(e^n)}V$ and 
with a nontrivial submodule of $\SL(V)$-invariants. Let $u$ be the smallest 
number for which the factor module $W'_{t+u}/W'_t$ has a nonzero module
 of $SL(V)$-invariants. Then the factor 
 $W'_{t+u}/W'_{t+u-1}=L_{\lambda (t+u)}$ contains a nonzero $SL(V)$-invariant.
  This is impossible, because the unique
  irreducible submodule of $L_{\lambda (t+u)}$
has  highest weight $\lambda (t+u)$ and is not one-dimensional.\end{proof}
  
\begin{corollary}
 Let $W$ be a polynomial homogeneous $\GL(V)$-module of degree $d=en$ with a good filtration. 
 Then the dimension of the $W^{\SL(V)}$-invariants in $W$ is equal to
  $n_{(e^n )}(W)$.
\end{corollary}

We continue with the application to quiver representations. Assume that
$Q$ is a quiver without oriented cycles.
The above reasoning generalizes directly to products of general linear groups, 
as the Schur modules for products of general linear groups are just 
tensor products of the Schur modules for the factors. In dealing below 
with the coordinate rings $K[\Rep(Q,\alpha )]$ and their good filtrations, 
notice we are really dealing with their homogeneous components which are 
finite dimensional representations.

\begin{proposition}\label{propGF3}
Let $Q$ be a quiver with no oriented cycles and let $\alpha$ be a dimension vector. The coordinate ring
$K[\Rep(Q,\alpha )]$ has a good filtration as a 
$\GL(Q,\alpha )$-module. Thus the dimension of the spaces of semi-invariants 
$\SI(Q,\alpha )_\sigma$ can be calculated as a multiplicity of the 
corresponding tensor products of Schur functors in the coordinate ring 
$K[\Rep(Q,\alpha )]$. In particular this dimension does not depend on the 
characteristic of $K$.
\end{proposition}
\begin{proof}
 The coordinate ring $K[\Rep(Q,\alpha )]$ has the following decomposition.
$$K[\Rep (Q,\alpha )]=\otimes_{a\in Q_1} \Sym (V(ta )\otimes V(ha)^*) .$$
Now, using the straightening law (comp. \cite{DC-E-P} or \cite[Theorem
(2.3.2)]{W}) we have that $\Sym (V(ta )\otimes V(ha)^*)$ has a characteristic 
free filtration with 
associated graded object 
$\oplus_{\lambda (a)} L_{\lambda (a)}V(ta)\otimes L_{\lambda (a)}V(ha)^*$. 
This is a good filtration. Applying Proposition~\ref{propGF1} 
we get that $K[\Rep(Q,\alpha )]$ has a good filtration as a 
$\GL(Q,\alpha )$-module. Now Proposition~\ref{propGF2} 
(for a product of general linear groups) gives the result.
\end{proof}

\section{A generalization to Covariants}\label{sec5}
The statement of Theorem~\ref{theo1} generalizes from semi-invariants to 
covariants.
Let us assume that $\langle \beta ,\gamma \rangle\ge
0$. Assume that the cycle $[q^{-1}(V)]$ for generic $V\in \Rep(Q,\alpha )$
decomposes as follows.
\begin{equation}\label{eqb}
[q^{-1}(V)]=\sum_{{\underline\lambda}:Q_0\rightarrow{\mathcal P}} 
N(\beta ,\alpha,{\underline\lambda})\prod_{x\in Q_0}
[{\overline{\underline\lambda}}(x)]_x.
\end{equation}
Here $\overline{\underline{\lambda}}(x)$ is the complement
of $\underline{\lambda}(x)$ in a $\beta(x)\times \gamma(x)$ rectangle.

Let $W(x)$ be a $\gamma(x)$-dimensional $K$-vector space for all $x\in Q_0$.
We can identify $\Rep(Q,\gamma)$ with $\bigoplus_{a\in Q_0}\Hom(W(ta),W(ha))$
and $\GL(Q,\gamma)$ with $\prod_{x\in Q_0} \GL(W(x))$.

We define $M(\beta,\alpha,\underline{\lambda})$ as the
multiplicity of $\det^{\sigma}$ in 
$$K[\Rep(Q,\gamma)]\otimes \bigotimes_{x\in Q_0}{\textstyle
\bigwedge^{\underline{\mu}(x)}}W(x).$$

\begin{proposition}\label{prop3}
Assume that $\sum_{x\in Q_0}
|{\underline\mu}(x)|=\langle \beta ,\gamma \rangle$. Then
$$N(\beta,\alpha, {\underline\mu})= M(\beta,\alpha,\underline{\mu}).$$

\end{proposition}

%
%
We will reduce Proposition~\ref{prop3} to Theorem~\ref{theo1}.
Let us write
$$\overline{\underline{\mu}} (x)= (\gamma(x)^{b_1(x)},(\gamma(x)-1)^{b_2(x)},\dots,
1^{b_{\gamma(x)}(x)})
$$
for all $x$, where $\overline{\underline{\mu}}(x)$ is the complement
of $\underline{\mu}(x)$ inside a $\beta(x)\times \gamma(x)$ rectangle.
We introduce the quiver $\widehat{Q}$ with
$$
\widehat{Q}_0=Q_0\cup \bigcup_{x\in
Q_0}\{y_{1,x},y_{2,x},\dots,y_{\gamma(x),x}\}.
$$
and
$$
\widehat{Q}_1=Q_1\cup \bigcup_{x\in Q_1}\{a_{1,x},a_{2,x},
\dots,a_{\gamma(x),x}\}
$$
where
$$
a_{i,x}:y_{i-1,x}\to y_{i,x}.
$$
for all $i$ and $x$. We use the convention $y_{0,x}=x$.

We define dimension vectors $\widehat{\beta},\widehat{\gamma}$
by
$$
\widehat{\beta}(x)=\beta(x),\quad \text{$x\in Q_0$,}
$$
$$
\widehat{\beta}(y_{i,x})=b_1(x)+b_2(x)+\cdots+b_{\gamma(x)-i+1},\quad
\text{$i=1,2,\dots,\gamma(x)$}
$$
and
$$
\widehat{\gamma}(x)=\gamma(x),\quad\text{$x\in Q_0$,}
$$
$$
\widehat{\gamma}(y_{i,x})=\gamma(x)-i+1, \quad\text{$i=1,2,\dots,\gamma(x)$}.
$$
\begin{lemma}
We have
$$
N(\beta,\alpha,\underline{\mu})=N(\widehat{\beta},\widehat{\alpha}).
$$
\end{lemma}
\begin{proof}
From (\ref{eq1}) follows that
$$
N(\widehat{\beta},\widehat{\alpha})=
\sum_{\underline\lambda :\widehat{Q}_1\rightarrow {\mathcal P}} \prod_{x\in
\widehat{Q}_0} \prod_{{\scriptstyle a\in \widehat{Q}_1\atop \scriptstyle ta=x}}
[{\underline\lambda} (a)]_x  \prod_{{\scriptstyle a\in \widehat{Q}_1\atop
\scriptstyle ha=x}}[{\overline{\underline\lambda }(a)}]_x.
$$
To get the class of a point at vertex $y_{\gamma(x),x}$ we must
have 
$$\underline{\lambda}(a_{\gamma(x),x})=(\gamma(y_{\gamma(x),x})^{\beta(y_{\gamma(x),x})})=(1^{b_1}).
$$
Now $\underline{\overline{\lambda}}(a_{\gamma(x),x})$ and 
$\underline{\lambda}(a_{\gamma(x),x})$ fit together into
a $\beta(y_{\gamma(x)-1,x})\times \gamma(y_{\gamma(x),x})=(b_1+b_2)\times 1$
rectangle. It follows that
$$
\underline{\overline{\lambda}}(a_{\gamma(x),x})=(1^{b_2}).
$$
To get a nonzero summand, $\underline{\overline{\lambda}}(a_{\gamma(x),x})$
and $\underline{\lambda}(a_{\gamma(x)-1,x})$ have to fit together
into a $\beta(y_{\gamma(x)-1,x})\times \gamma(y_{\gamma(x)-1,x})=
 (b_1+b_2)\times 2$ rectangle. We see that
$$
\underline{\lambda}(a_{\gamma(x)-1,x})=(2^{b_1},1^{b_2}).
$$
Continuing by induction, we see that
$$
\underline{\lambda}(a_{1,x})=(\gamma(x)^{b_1(x)},\dots,1^{b_{\gamma(x)}(x)})=
\overline{\underline\mu}(x).
$$
We have
$$
N(\widehat{\beta},\widehat{\alpha})=
\prod_{a\in \widehat{Q}_1} \big(\sum_{\lambda}
[\lambda]_{ta}[{\overline\lambda}]_{ha} \big)=$$
$$
=
\big(\prod_{a\in Q_1} \big(\sum_{\lambda}
[\lambda]_{ta}[{\overline\lambda}]_{ha} \big)\big)
\big(\prod_{a\in \widehat{Q}_1\setminus Q_1} \big(\sum_{\lambda}
[\lambda]_{ta}[{\overline\lambda}]_{ha} \big)\big)
$$
Using our calculations for $\underline{\lambda}(a_{i,x})$ above,
we see that this is equal to
$$
\big(\prod_{a\in Q_1} \big(\sum_{\lambda}
[\lambda]_{ta}[{\overline\lambda}]_{ha} \big)\big)\prod_{x\in
Q_0}[\underline{\mu}(x)]_x.
$$
From (\ref{eq0}) and (\ref{eqb}) follows that
$$
N(\widehat{\beta},\widehat{\alpha})=
\big(\sum_{{\underline\lambda}:Q_0\rightarrow{\mathcal P}} N(\beta ,\alpha
, {\underline\lambda})\prod_{x\in Q_0} [\overline{\underline\lambda}(x)]_x\big)
\prod_{x\in
Q_0}[\underline{\mu}(x)]_x=N(\beta,\alpha,\underline{\mu}).
$$
\end{proof}
\begin{lemma}
We have
$$
M(\beta,\alpha,\underline{\mu})=M(\widehat{\beta},\widehat{\alpha}).
$$
\end{lemma}
\begin{proof}
The proof goes similar to the proof of the previous lemma. From
(\ref{eqstarstar2}) follows that
$$
M(\widehat{\beta},\widehat{\alpha})=\sum_{\underline\lambda:\widehat{Q}_1\to
{\mathcal P}} \prod_{x\in
\widehat{Q}_0} \mult\big({\textstyle \bigwedge^{(\gamma (x)^{\beta (x)})}}W(x) 
; (
\bigotimes_{\scriptstyle a\in \widehat{Q}_1\atop \scriptstyle ta =x}
{\textstyle \bigwedge^{{\underline\lambda}(a)}}W(x))\otimes (
\bigotimes_{\scriptstyle a\in \widehat{Q}_1\atop \scriptstyle ha=x}
{\textstyle \bigwedge^{{\overline{\underline\lambda}(a)}}}W(x))\big ).
$$
To get a nonzero summand, we get the same conditions for
$\underline{\lambda}(a_{i,x})$ as in the previous lemma. We obtain
$$
M(\widehat{\beta},\widehat{\alpha})=
$$
$$=
\sum_{\underline\lambda:Q_1\to
{\mathcal P}} \prod_{x\in
Q_0} \mult\Big({\textstyle \det_x^{\beta(x)}} ; \big(
\bigotimes_{\scriptstyle a\in Q_1\atop \scriptstyle ta =x}
{\textstyle \bigwedge^{{\underline\lambda}(a)}}W(x)\big)\otimes \big(
\bigotimes_{\scriptstyle a\in Q_1\atop \scriptstyle ha=x}
{\textstyle \bigwedge^{{\overline{\underline\lambda}(a)}}}W(x)\big)\otimes
{\textstyle \bigwedge^{\underline{\mu}(x)}}W(x)\Big)=
$$
$$
=\sum_{\underline\lambda:Q_1\to
{\mathcal P}} \prod_{x\in
Q_0} \mult\Big({\textstyle\det_x^{\beta(x)}} ; \big(
\bigotimes_{\scriptstyle a\in Q_1\atop \scriptstyle ta =x}
{\textstyle \bigwedge^{{\underline\lambda}(a)}}W(x)\big)\otimes \big(
\bigotimes_{\scriptstyle a\in Q_1\atop \scriptstyle ha=x}
{\textstyle \bigwedge^{\underline{\lambda}(a)}}W^\star(x)\otimes
{\textstyle\det_x^{\beta(ta)}}\big)\otimes \textstyle\bigwedge^{\underline{\mu}(x)}W(x)\Big)=
$$
$$
=\mult\Big({\textstyle \det^\sigma};
\bigoplus_{\underline{\lambda}:Q_1\to {\mathcal P}}
 \bigotimes_{a\in Q_1}
\big({\textstyle \bigwedge^{{\underline\lambda}(a)}}W(ta)\otimes 
{\textstyle\bigwedge^{{\underline{\lambda}(a)}}}W^\star(ha)\big)
\otimes \bigotimes_{x\in Q_0}{\textstyle
\bigwedge^{\underline{\mu}(x)}}W(x)\Big)=
$$
$$
=\mult\Big({\textstyle\det^{\sigma}};\Rep(Q,\gamma)\otimes\bigotimes_{x\in Q_0} {\textstyle\bigwedge^{\underline{\mu}(x)}}W(x)\Big)=M(\beta,\alpha,\underline{\mu}).
$$
\end{proof}

\section{Applications}\label{sec6}

Theorem~\ref{theo1} also allows us to exhibit an explicit basis of the
weight space $\SI(Q,\gamma )_{\langle\beta ,\cdot\rangle}$.

\begin{corollary}\label{cor1} 
Let $Q,\alpha ,\beta ,\gamma$ be as in Theorem~\ref{theo1}.
Assume that the general representation $V$ of dimension
$\alpha$ has $k$ subrepresentations of dimension
$\beta$. There exists a
nonempty Zariski open set $U$ in $\Rep_K(Q,\alpha )$ such that for $V\in
U$ the  semi-invariants $c^{V_1},\ldots ,c^{V_k}$ form a basis
in $\SI(Q,\gamma )_{\langle\beta ,\cdot\rangle}$, where $V_1,V_2,\dots,V_k$
are the subrepresentations of $V$ of dimension $\beta$.
\end{corollary}

\begin{proof}  Let us choose $V\in \Rep_K(Q,\alpha )$ such that
$q^{-1}(V)$ consists of $k$ points. Let $V_1 ,\ldots ,V_k$ be
the corresponding subrepresentations of $V$ of dimension $\beta$. It is
enough to prove that
$c^{V_1},\ldots ,c^{V_k}$ are linearly independent in $\SI(Q,\gamma
)_{\langle\beta ,\cdot\rangle}$. Let us consider the exact sequences
$$0\rightarrow V_i\rightarrow V\rightarrow W_i\rightarrow 0$$
for $i=1,\ldots ,k$. It is clear that for $i\ne j$ we have $\Hom_Q(V_i ,W_j
)\ne 0$. Therefore $c^{V_i}(W_j )=0$ for $i\ne j$. Therefore
it is enough to show that $\Hom_Q (V_i ,W_i )=0$ for $i=1,\ldots ,k$. In
\cite{S1} it was proved that $\Hom_Q(V_i ,W_i )$ is the tangent
space to the fiber of the map $q:Z(Q,\beta,\alpha)\rightarrow
\Rep_K(Q,\alpha )$ at the point $(V, V_i )$. The map $q$ is dominant
and generically it is $k:1$. Moreover, it is shown in \cite{CB} that the
differential $Dq$ is generically surjective. Therefore it
is generically an isomorphism because $Z(Q,\beta,\alpha)$ is smooth.
Therefore for $V$ in
some nonempty Zariski open set we have $\Hom_Q (V_i ,W_i )=0$ for $i=1,\ldots ,k$.
\end{proof}
\begin{proposition}
Let $Q,\alpha ,\beta ,\gamma$ be as in Theorem~\ref{theo1}.
Let $U_\beta\subseteq \Rep(Q,\beta)$ and $U_\gamma\subseteq\Rep(Q,\gamma)$
be nonempty Zariski open subsets,
stable under $\GL(\beta)$ and $\GL(\gamma)$
respectively. Then there exists a nonempty
Zariski open subset $U$ of $\Rep(Q,\alpha)$ such that for all $V\in U$
we have that every 
$\beta$-dimensional subrepresentation of $V$ lies in $U_\beta$
and every $\gamma$-dimensional factor representation of $V$ lies in $U_\gamma$.
\end{proposition}
\begin{proof}
Let 
$$D\subset Z(Q,\beta,\alpha)\subset \Rep(Q,\alpha)\times
\Grass(\beta,\alpha)$$ 
be the subset of pairs $(V,W)$ such that $W$ does not
lie in $U_{\beta}$ or $V/W$ does not lie in $U_{\gamma}$.
We claim that $D$ is a Zariski closed
strict subset of $Z(Q,\beta,\alpha)$.
For every $z\in \Grass(\beta,\alpha)$ we can choose
local sections $e_1,e_2,\dots,e_{\alpha(x)}$ and an
open neighborhood $X$ of $z$ such
that
$e_1(W),\dots,e_{\beta(x)}(W)$
are a basis of $W$ and
$e_1(W),\dots,e_{\alpha(x)}(W)$ are a basis of $K^{\alpha(x)}$
for all $W\in X$.
With respect to this basis, $V$ has the form
$$
V(a)=
\begin{pmatrix}
V'(a) & \star\\
0 & V''(a)\end{pmatrix}
$$

where $V'\in \Rep(Q,\beta)$, $V''\in \Rep(Q,\gamma)$ and $\star$ is an arbitrary matrix.

We define
$$
r_X:p^{-1}(X)\to \Rep(Q,\beta)\times\Rep(Q,\gamma)\times X.
$$
by
$$
(V,W)\mapsto (V',V'',W),
$$
It is clear that
$$
D\cap p^{-1}(X)=r_X^{-1}((D_\beta\times \Rep(Q,\gamma)\times X)\cup
(\Rep(Q,\beta)\times D_{\gamma}\times X))
$$
where $D_\beta$ and $D_\gamma$ are the complements 
of $U_\beta$ and $U_\gamma$ respectively.
Therefore we have that $D\cap p^{-1}(X)$ is a Zariski closed proper subset of 
$p^{-1}(X)$.
Since such $p^{-1}(X)$ cover $Z(Q,\beta,\alpha)$, we conclude
that $D$ is Zariski closed proper subset of $Z(Q,\beta,\alpha)$.
For some $U'\subset \Rep(Q,\alpha)$ open nonempty, $q^{-1}(V)$ is finite
for all $V\in U'$. Take $U=U'\setminus q(D)$. The map $q$ is proper
so $q(D)$ is closed and $U$ is therefore open. We also claim that
$U$ is nonempty. Indeed, the restriction $q:q^{-1}(U')\to U'$
is quasi-finite, so $q(q^{-1}(U')\cap D)\subseteq U'\cap q(D)\neq U'$
because the dimension of $q^{-1}(U')\cap D$ is strictly
smaller than $\dim q^{-1}(U')=\dim U'$.
If $V\in U$ then $q^{-1}(V)$ is finite and
for every $(V,W)\in q^{-1}(V)$ we have that $W$ and $V/W$ 
lie in $U_\beta$ and $U_{\gamma}$ respectively.
\end{proof}
The proposition can be roughly reformulated as follows. If
$V$ is a general representation of dimension $\alpha=\beta+\gamma$
with $\langle\beta,\gamma\rangle =0$, then all subrepresentations
of dimension $\beta$ and all factor representation of dimension $\gamma$
are in general position as well. 
\begin{corollary}
Suppose that $\beta,\gamma,\delta$ are dimension vectors
such that $\langle
\beta,\gamma\rangle=\langle\beta,\delta\rangle=\langle\gamma,\delta\rangle=0$.
Then we have the following equality
$$
N(\beta,\beta+\gamma)N(\beta+\gamma,\beta+\gamma+\delta)=
N(\beta,\beta+\gamma+\delta)N(\gamma,\gamma+\delta).
$$
\end{corollary}
\begin{proof}
With the previous proposition this is now a simple counting argument.
Choose $V\in \Rep(Q,\beta+\gamma+\delta)$ in general position.
We count the number of pairs $(V_1,V_2)$ such that $V_1$
is an $\beta$-dimensional subrepresentation of $V_2$
and $V_2$ is a $(\beta+\gamma)$-dimensional subrepresentation of $V$.
On the one hand, $V$ has $N(\beta+\gamma,\beta+\gamma+\delta)$
$(\beta+\gamma)$-dimensional subrepresentations $V_2$ and each such 
subrepresentation
(since it is again in general position) has exactly
$N(\beta,\beta+\gamma)$ subrepresentations $V_1$ of dimension $\beta$.

On the other hand $V$ has $N(\beta,\beta+\gamma+\delta)$
$\beta$-dimensional
subrepresentations $V_1$. For each $V_1$, $V/V_1$ is again
in general position and $V/V_1$ has exactly $N(\gamma,\gamma+\delta)$
subrepresentations of dimension $\gamma$. Note also that there is a $1$--$1$ correspondence
between $\gamma$-dimensional subrepresentations of $V/V_1$
and $\beta+\gamma$ dimensional subrepresentations $V_2$ of $V$ containing
$V_1$.
Comparison of the two computations completes the proof.
\end{proof}

\begin{example}
 Let $Q=\theta (m)$ be the quiver with two vertices $x,
y$ and $m$ arrows $a_1 ,\ldots ,a_m$, with $ta_i =x, ha_i =y$ for $i=1,\ldots ,m$.
Assume $m=2r$ is even.
Let $\beta$ be the dimension vector $\alpha (x)=\alpha (y)=r+1$. Consider
the subdimension vector $\beta $ with $\beta (x)=1, \beta
(y)=r$. Then $\langle\beta ,\alpha-\beta\rangle =0$. The
Littlewood-Richardson calculation shows that $N(\beta ,\alpha )= {2r \choose
r}$.

In particular, for $r=2$, we have $Q= \theta (4)$ is the quiver with two
vertices $x, y$ and four arrows $a, b, c, d$ from $x$ to $y$.
Consider the dimension vector $\alpha$ with $\alpha (x)=\alpha (y) =3$. Consider the map
$$C: \Rep(Q, \alpha)\rightarrow \A^{20}$$
where we identify $\A^{20}$ with the space of
quaternary cubics.
The map $C$ sends a representation $V$ to the point $\det\ (X_0  V(a)+X_1
V(b)+X_2 V(c)+X_3 V(d))$.

Consider the subrepresentations of dimension $\beta$ where $\beta (x)=1,
\beta (y)=2$.
Our calculation tells us that a generic representation of dimension
$\alpha$ has $6$ subrepresentations of
dimension $\beta$. They correspond to six lines on the cubic surface 
defined by $C(V)$. These lines
 represent  the cubic surface as a projective plane with $6$ points 
 blown
up. Given a subrepresentation $W$ of dimension $\beta$, the corresponding
line is constructed as follows. Choose bases of $V(x)$ and $V(y)$ such
that $W(x)$ is spanned by the first basis vector and $W(y)$ is
spanned by the first 2 basis vectors. Then the cubic surface
is defined by 
$$
\det\ \begin{pmatrix}l_{1.1}&l_{1,2}&l_{1,3}\\
l_{2,1}&l_{2,2}&l_{2,3}\\
0&l_{3,2}&l_{3,3}\end{pmatrix}=0,$$
where the $l_{i,j}$ are linear functions.
Now $l_{1,1}=l_{2,1}=0$ defines a line on the surface.
\end{example}

\end{document}